\documentclass[10pt]{article}

\usepackage[english]{babel}

\usepackage{algpseudocode}
\usepackage{amsmath}
\usepackage{amsthm}
\usepackage{amsfonts}
\usepackage{amssymb}
\usepackage{color}
\usepackage{textcomp}

\usepackage{url}

\usepackage{geometry}

\theoremstyle{plain}
\newtheorem{thm}{Theorem}
\newtheorem{lem}[thm]{Lemma}
\newtheorem{prop}[thm]{Proposition}
\newtheorem{cor}{Corollary}
\theoremstyle{definition}
\newtheorem{defn}{Definition}

\newtheorem{exmp}{Example}
\newtheorem{open}{Open problem}
\newtheorem{prob}{Problem}
\theoremstyle{remark} 
\newtheorem{rem}{Remark}

\allowdisplaybreaks

\title{An Equivalent Representation of Generalized Differentials}
\author{Valentin Suder}
\date{\today}

\begin{document}

\maketitle

\begin{abstract} 
    We propose an equivalent formula for the higher-order derivatives used in the study of Generalized Almost Perfect Nonlinear functions over an arbitrary finite field of characteristic $p$. The result is obtained by counting the number of subsets of the prime field with a fixed cardinality for which the sum of their elements is constant. We then ask related questions regarding the diversity of higher-order derivatives.
\end{abstract}

\section{Introduction}
Almost Perfect Nonlinear (APN) functions were brought to light by Nyberg~\cite{Nyberg93} for granting the symmetric cryptosystems they compose with the best possible resistance to differential cryptanalysis~\cite{BihamS90}. APN functions, and more generally Perfect Nonlinear functions,  are mathematical objects with extremal properties and are therefore attracting a continuous research effort from different mathematical communities (see~\cite{CCZ98,dempwolff2014dimensional,PUB16} for some examples).

In 2017, Kuroda and Tsujie~\cite{KurodaT17} introduced the concept of Generalized Almost Perfect Nonlinear (GAPN) functions over an arbitrary finite field. They argue that APN functions over a binary field are intuitively closer as a notion to GAPN functions than they are to APN functions over a finite field with odd characteristic where Perfect Nonlinear functions exist. 

Later, \"Ozbudak and Salagean~\cite{OzbudakS21,SalageanO23}, reinterpreted this generalization by mean of discrete derivatives, and more specifically higher-order derivatives of a very particular forms. They were therefore successful in finding new classes of GAPN functions using similar tooling as for the search of APN functions over binary fields.

In this article, we give a different formula from \"Ozbudak and Salagean, by using different higher-order discrete derivatives. More formally, we give the proof of the following Lemma.
\begin{lem}\label{thm:main}
    Let $F:\mathbb{F}_{p^n}\rightarrow \mathbb{F}_{p^n}$, and $\alpha\in\mathbb{F}_{p^n}$, then for all $x\in\mathbb{F}_{p^n}$ we have
    \[\Delta_{\underbrace{\alpha,\alpha,\dots,\alpha}_{p-1}} F(x)= -\Delta_{\alpha,2\alpha,\dots,(p-1)\alpha} F(x).\]
\end{lem}

Although the result itself is somewhat anecdotic in the sense that it does not seem to bring any new insight on generalized APN functions, the method to obtain the formula sparks some independent interest. Firstly, it relies on identifying the cardinality of the sets:
\begin{align*} 
S_\ell(u) &= \{S\subseteq \mathbb{F}_p \mid \#S=\ell,~\sum_{s\in S} s = u\}, \quad \mbox{and} \\
S_\ell^*(u) &= \{S\subseteq \mathbb{F}_p^* \mid \#S=\ell,~\sum_{s\in S} s = u\}
\end{align*}
for any $u\in\mathbb{F}_p$ and any $0\leq \ell \leq p$ (resp. $0\leq \ell < p$). While the cardinality of $S_\ell(u)$ can be found easily, the cardinality of $S_\ell^*(u)$ seems to be less trivial to obtain.

Secondly, the equation from Lemma~\ref{thm:main} raises a more general question about the diversity of higher-order derivatives for a given function, and up to Extended-Affine  (EA) equivalence. For a given function, all of its higher-order derivatives with an order higher than the function's algebraic degree trivially vanish (given that the algebraic degree of the first-order derivative is at least one less than the degree of the function being differentiated). Lemma~\ref{thm:main} hints that different higher-order derivatives, i.e. in different directions, can be functionally equal and non-zero.

This article is organized as follows: Section~\ref{sec:prelim} introduces the formal definitions and preparatory results, Section~\ref{sec:lset} gives a treatment of the challenges to achieve the main result. Finally, before the concluding remarks, we propose some questions and answers regarding matching higher-order derivatives of distinct functions.

\section{Preliminaries}\label{sec:prelim}
In the remaining of this work, we denote by $\mathbb{F}_q$ the finite fields with $q$ elements, a power of a prime $p$ called the characteristic. The non-zero elements of a set $S$ forms a set denoted by $S^*$.
For a function $F:\mathbb{F}_q \rightarrow \mathbb{F}_q$ and an element $\alpha\in\mathbb{F}_q$, the discrete derivative of $F$ in the direction $\alpha$ is usually defined as:
\[\Delta_\alpha F(x) = F(x+\alpha) - F(x).\]

Functions for which the image set of all of their derivatives are as large as possible are called Perfect Nonlinear (PN) functions, that is when all of their derivatives are bijective functions of $\mathbb{F}_{p^n}$. When the characteristic of the field $p=2$, the relation $\Delta_\alpha F(x)=\Delta_\alpha F(x+\alpha)$ indicates that the largest image set a derivative can have is the size of an hyperplane, i.e. $2^{n-1}$ and functions with such derivatives are called Almost Perfect Nonlinear (APN) functions. 
Hence PN functions only exist when $p$ is odd. Note that when $p$ is odd, APN functions are forcefully defined to be the functions for which every element of the field has at most 2 pre-images by each of their derivatives. 

In 1994, Lai~\cite{Lai94} extended the concept of discrete derivatives. The $d$-order derivative of $F$ along the directions $\alpha_i\in\mathbb{F}_q$, where $i=1,\dots,d$ is denoted $\Delta_{\alpha_1,\dots,\alpha_d}F$ and satisfies:
\begin{equation}\label{eq:higher-order}
\Delta_{\alpha_1,\dots,\alpha_d} F(x) = \Delta_{\alpha_1}\Delta_{\alpha_2,\dots,\alpha_d} F(x) = \sum_{j=0}^{d} (-1)^{d-j} \sum_{\{i_1,\dots,i_j\}\subseteq \{1,\dots,d\}} F(x + \sum_{k=1}^j \alpha_{i_k}).
\end{equation}

In the following proposition, we list some noticeable properties concerning these derivatives, and we leave the proofs as an exercise.

\begin{prop}
    Let $F$ be a function over $\mathbb{F}_{p^n}$. Let $\alpha,\beta\in\mathbb{F}_{p^n}$ and let $\iota$ be an element of the base field $\mathbb{F}_p$. Then, for all $x\in\mathbb{F}_{p^n}$ we have:
    \begin{enumerate}
        \item $\Delta_{\alpha,\beta} F(x) = \Delta_{\beta,\alpha}F(x)$
        \item $\Delta_\alpha F(x) = - \Delta_{-\alpha} F(x+\alpha)$
        \item $\Delta_\alpha F(x) - \Delta_\beta F(x) = \Delta_{\alpha-\beta} F(x+\beta)$
        \item $\Delta_{\iota\alpha} F(x) = \Delta_\alpha \sum_{j=0}^{\iota-1} F(x+j\alpha)$
    \end{enumerate}
\end{prop}

In~\cite{KurodaT17}, Kuroda and Tsujie redefined the concept of Almost Perfect Nonlinear functions for finite fields of odd characteristic. They consider that a function $F:\mathbb{F}_{p^n}\rightarrow\mathbb{F}_{p^n}$ should be APN when for all $\alpha\neq0,\beta\in\mathbb{F}_{p^n}$ the following equation:
\begin{equation}
    \sum_{\iota\in\mathbb{F}_p} F(x +\iota\alpha) = \beta
\end{equation}
has at most $p$ solutions. In this case, the function $F$ is said to be a Generalized APN (APN) function as it coincides with the concept of APN function when $p=2$ and generalized intuitively the same concept for $p$ odd.

Later, \"Ozbudak and Salagean~\cite{OzbudakS21} noticed that Generalized APN functions can be worked out using a specific type of higher-order differential, that they conveniently call generalized differentials.
\begin{lem}[\cite{OzbudakS21} Lemma 1]\label{lem:integral}
    For any function $F:\mathbb{F}_{p^n}\rightarrow\mathbb{F}_{p^n}$ and any $\alpha\in\mathbb{F}_{p^n}$ we have
    \[\Delta_{\underbrace{\alpha,\dots,\alpha}_{p-1}}F(x) = \sum_{\iota\in\mathbb{F}_p} F(x+\iota \alpha).\]
\end{lem}

\section{Proof of the main result}\label{sec:lset}
We are now establishing the proof of Lemma~\ref{thm:main}. The strategy we employ is to use Lemma~\ref{lem:integral} in order to show that: 
\begin{equation}
    \Delta_{\alpha, \dots, (p-1)\alpha} F(x) = - \sum_{\iota\in\mathbb{F}_p} F(x+\iota\alpha).
\end{equation}

By using Equation~\eqref{eq:higher-order}, for any $\alpha\in\mathbb{F}_{p^n}$ we have that
\begin{equation}\label{eq:spatial}
    \Delta_{\alpha, \dots, (p-1)\alpha} F(x) = \sum_{\ell=0}^{p-1} (-1)^{p-1-\ell} \sum_{\{k_1,\dots,k_\ell\} \subseteq \mathbb{F}_p^*} F(x + \alpha\sum_{i=1}^\ell k_i).
\end{equation}

Note that the term $\sum_{i=1}^\ell k_i$ can only be one of the $p$ elements in $\mathbb{F}_{p^n}$ for any $\ell$.  Note also that Equation~\eqref{eq:spatial} is an alternated sum running through every subsets of $\mathbb{F}_p^*$ exactly once.
Our first step is then to count, for a fixed $0\leq \ell<p$, how many times the sum $\sum_{i=1}^\ell k_i$, $k_i\in\mathbb{F}_p^*$, is equal to a certain $u\in\mathbb{F}_{p^n}$. In other words, we need to find the cardinality of the previously defined set:
\[ S_\ell^*(u) = \{S\subseteq \mathbb{F}_p^* \mid \#S=\ell,~\sum_{s\in S} s = u\},\]
through the cardinality of its closely related set:
\[S_\ell(u) = \{S\subseteq \mathbb{F}_p \mid \#S=\ell,~\sum_{s\in S} s = u\}.\]

We describe these preliminary results in the next subsection.
Once the cardinality of the aforementioned sets are established, we will come back to the main proof in the Subsection~\ref{sec:main}.

 \subsection{Computing the cardinality of $S_\ell^*(u)$}
It is quite clear that every subset of $\mathbb{F}_p$ (resp. $\mathbb{F}_p^*$) with cardinality $0\leq \ell<p$ belongs to exactly one of the sets $S_\ell(u)$ (resp. $S_\ell^*(u)$). Hence, we deduce the following partitions:
\begin{align*}
    \bigcup_{u\in\mathbb{F}_p} S_\ell(u) &= \{ S\subseteq \mathbb{F}_p \mid \# S = \ell\},\\
    \bigcup_{u\in\mathbb{F}_p} S_\ell^*(u) &= \{ S\subseteq \mathbb{F}_p^* \mid \# S = \ell\},
\end{align*}
and their cardinality:
\begin{align}
    \# \bigcup_{u\in\mathbb{F}_p} S_\ell(u) &= \binom{p}{\ell} = \sum_{u\in\mathbb{F}_p} \# S_\ell(u),\\
    \# \bigcup_{u\in\mathbb{F}_p} S_\ell^*(u) &= \binom{p-1}{\ell} = \sum_{u\in\mathbb{F}_p} \# S_\ell^*(u). \label{eq:partition*}
\end{align}

\begin{prop}
    Let $0\leq \ell<p$. The sets $S_\ell(u)$ have the same cardinality for all $u\in\mathbb{F}_p$:
    \[ \# S_\ell(u) = \binom{p}{\ell}/p, \quad \forall u \in \mathbb{F}_p.\]
\end{prop}
\begin{proof}
    Let $S\in S_\ell(u)$, then for any $i=0,\dots,p-1$ we have that
    \[ i+S := \{i +s \mid \forall s\in S\} \in S_\ell(u+i\ell).\]
    Moreover, because $\ell$ and $p$ are necessarily prime to each other, for any $u\in\mathbb{F}_p$ and for any distinct $S,S'\in S_\ell(u)$, the sets $i+S,\, i+S' \in S_\ell(u+i\ell)$ are distinct. By double inclusion, we conclude that the sets $S_\ell(u)$ and $S_\ell(u+i\ell)$ have the same cardinality.
\end{proof}

The same arguments cannot be reused because none of the sets in $S_\ell^*(u)$ contain 0. Indeed, there will always  be some of the translated sets $i+S^*$ for some $i\in\mathbb{F}_p^*$ and $S^*\in S_\ell^*(u)$ containing the zero element. However, because we now deal with sums over the multiplicative group $\mathbb{F}_p^*$, a different approach can be taken. By using a recurrence relationships over $\ell$, we  show that only the cardinality of $S_\ell^*(0)$ differs from $\# S_\ell^*(u)$, $u\neq 0$.

\begin{prop}\label{prop:sl*}
    Let $0\leq \ell < p$. The sets $S_\ell^*(u)$ have the same cardinality for all $u\in\mathbb{F}_p^*$:
    \[S_\ell^* := \# S_\ell^*(u), \quad \forall u \in \mathbb{F}_p^*.\]
\end{prop}
\begin{proof}
    Since there is no divisors of zero in $\mathbb{F}_p^*$, we can be certain that for all $S\in S_\ell^*(u)$, $u\neq 0$, and for any $t\in\mathbb{F}_p^*$:
    \[ t\times S := \{ts \mid s\in S\} \in S_\ell^*(tu).\] 
    We conclude by double inclusion.
\end{proof}

In the following we also use the notation $S_\ell^{*0}:=\#S_\ell^*(0)$ for practicality.
By using the result and notation from the proposition above, we can re-write Equation~\eqref{eq:partition*}:
\[\binom{p-1}{\ell} = S_\ell^{*0} + (p-1)S_\ell^*.\]
In this next propositions, we reuse the same notations as in Proposition~\ref{prop:sl*}.
\begin{prop}\label{prop:rec_ell}
    For any $u\in\mathbb{F}_p$, and $0\leq \ell <p$, we have 
    \[\# S_\ell(u) = \binom{p}{\ell}/p = S_\ell^* + S_{\ell-1}^* = S_\ell^{*0} + S_{\ell-1}^{*0}.\]
\end{prop}
\begin{proof}
    Since there can only be zero or one appearance exactly of the element 0 in each of the sets $S\in S_\ell(u)$, we have:
    \begin{align*}
        S_\ell(u) &= \{S\subseteq \mathbb{F}_p \mid \#S=\ell, \sum_{s\in S} s=u\}\\
        &= \{S\subseteq \mathbb{F}_p^* \mid \#S=\ell, \sum_{s\in S} s=u\} \\
        & \qquad\bigcup~ \{ \{0\} \cup S \mid S\subseteq\mathbb{F}_p^*,~ \#S=\ell-1,~ \sum_{s\in S} s=u\} \\
        &= S_\ell^*(u)~ \bigcup~ \{ \{0\} \cup S \mid S\in S_{\ell-1}^*(u)\}.
    \end{align*}
\end{proof}

\begin{prop}\label{prop:diff-1}
    Let $0\leq \ell < p$. We have
    \[S_\ell^{*0} - S_\ell^* = (-1)^\ell.\]
\end{prop}
\begin{proof}
    First, we remark that $S_0^*(0) = \{\emptyset\}$ and $S_0^*(u)=\{\}$, for $u\neq0$, meaning that their respective cardinality are $S_0^{*0} =1$ and $S_0^* = 0$. Similarly, we have $S_1^{*0}=0$ and $S_1^*=1$. Now we use the recurrence relation from Proposition~\ref{prop:rec_ell} repetitively:
    \begin{align*}
        S_\ell^{*0} &= \binom{p}{\ell}/p - S_{\ell-1}^{*0} \\
        &= \frac{\binom{p}{\ell} - \binom{p}{\ell-1}}{p} + S_{\ell-2}^{*0} \\
        & \qquad \vdots \\
        &= \frac{\binom{p}{\ell} - \dots + (-1)^{\ell-1}\binom{p}{1}}{p} +(-1)^{\ell} S_0^{*0} \\
        &= \frac{\binom{p}{\ell} - \dots + (-1)^{\ell-1}\binom{p}{1}}{p} +(-1)^\ell \\
        &= \frac{1}{p} \sum_{i=0}^{\ell-1} (-1)^i \binom{p}{\ell-i} + (-1)^\ell. 
    \end{align*}
Similarly,
    \begin{align*}
        S_\ell^* &= \frac{\binom{p}{\ell} - \dots + (-1)^{\ell-1}\binom{p}{1}}{p} +(-1)^{\ell} S_0^{*} \\
        &= \frac{\binom{p}{\ell} - \dots + (-1)^{\ell-1}\binom{p}{1}}{p} \\
        &= \frac{1}{p} \sum_{i=0}^{\ell-1} (-1)^i \binom{p}{\ell-i}.
    \end{align*}
    We conclude this proof by taking the difference between the two relations above.
\end{proof}

Finally, the capstone of Lemma~\ref{thm:main}'s proof, that is concluded in the next subsection, resides in the following corollary. At this stage, the proof should be a walk in the park as it is purely calculation, but we still include it for clarity.
\begin{prop}\label{prop:congruence-1}
    With the above notations, we have:
    \[\sum_{i=1}^{(p-1)/2} S_{2i}^* - S_{2i-1}^* \equiv - 1\pmod{p}.\]
\end{prop}
\begin{proof}
    First, we note that $\binom{p}{\ell}/p = \binom{p-1}{\ell-1}/\ell$, for $0<\ell<p$: Proposition~\ref{prop:rec_ell} implies the following induction:
    \begin{align*}
  S_\ell^* - S_{\ell-1}^* &= \frac{\binom{p}{\ell}}{p} - 2S_{\ell-1}^* \\
  &= \frac{\binom{p-1}{\ell-1}}{\ell} - 2\left( \frac{\binom{p-1}{\ell-2}}{\ell-1} - S_{\ell-2}^*\right) \\
  &= \frac{\binom{p-1}{\ell-1}}{\ell} - 2\left( \frac{\binom{p-1}{\ell-2}}{\ell-1} - \frac{\binom{p-1}{\ell-3}}{\ell-2} + \dots - S_{1}^*\right) \\
  &\equiv \frac{(-1)^{\ell-1}}{\ell} - 2\frac{(-1)^{\ell-2}}{\ell-1} + 2\frac{(-1)^{\ell-3}}{\ell-2} - \dots + 2(-1)^{\ell-1-(\ell-1)} \pmod p\\
  &\equiv (-1)^{\ell-1}\times (\frac{1}{\ell} + \frac{2}{\ell-1} + \frac{2}{\ell-2} + \dots + 2) \pmod p\\
  &\equiv (-1)^{\ell-1} (\frac{1}{\ell} + \sum_{i=1}^{\ell-1} \frac{2}{\ell-i}) \pmod p
\end{align*}

Hence by summation,

\begin{align*}
  \sum_{i=1}^{(p-1)/2} S_{2i}^* - S_{2i-1}^* &\equiv \sum_{i=1}^{(p-1)/2} - \frac{1}{2i} - \sum_{j=1}^{2i-1} \frac{2}{2i-j} \pmod p\\
  &\equiv - (\frac{1}{p-1} + \frac{2}{p-2} + \frac{3}{p-3} + \dots + p-1) \\
  &\equiv - \sum_{i=1}^{p-1} \frac{i}{p-i} \equiv - \sum_{i=1}^{p-1} i(p-i)^{-1} \equiv - \sum_{i=1}^{p-1} -i\times i^{-1}\\
  &\equiv - \sum_{i=1}^{p-1} -1 \equiv -1 \pmod p.
\end{align*}
\end{proof}

\subsection{Proof of Lemma~\ref{thm:main}}\label{sec:main}
We are now using the results of the previous subsection in order to simplify Equation~\eqref{eq:spatial} into the result. We can first specify Equation~\eqref{eq:spatial} with $\alpha =1$ without loss of generality and then it becomes:
\begin{align*}
    \Delta_{1,\dots,p-1} F(x) &= \sum_{\ell=0}^{p-1} (-1)^{p-1-\ell} \sum_{\{k_1,\dots,k_\ell\}\subseteq \mathbb{F}_p^*} F(x+\sum_{i=0}^\ell k_i) \\
    &= \sum_{\ell=0}^{p-1} (-1)^\ell \left( S_\ell^{*0} F(x) + \sum_{u\in\mathbb{F}_p^*} S_\ell^* F(x+u)\right) \\
    &= S_0^{*0}F(x) + \sum_{i=1}^{(p-1)/2} \left( S_{2i}^{*0} F(x) + \sum_{u\in\mathbb{F}_p^*} S_{2i}^*F(x+u)\right) \\
    & \qquad - \sum_{i=1}^{(p-1)/2} \left( S_{2i-1}^{*0} F(x) + \sum_{u\in\mathbb{F}_p^*} S_{2i-1}^*F(x+u)\right) \\
    &= F(x)\left( 1+ \sum_{i=1}^{(p-1)/2} S_{2i}^{*0} - S_{2i-1}^{*0}\right) + \sum_{u\in\mathbb{F}_p^*} F(x+u)\sum_{i=1}^{(p-1)/2} S_{2i}^* - S_{2i-1}^*.
\end{align*}

Proposition~\ref{prop:diff-1} informs us that 
\begin{equation*}
    1+ \sum_{i=1}^{(p-1)/2} S_{2i}^{*0}- S_{2i-1}^{*0} = \sum_{i=1}^{(p-1)/2} S_{2i}^* - S_{2i-1}^*.
\end{equation*}

We conclude the proof of Lemma~\ref{thm:main} by using Proposition~\ref{prop:congruence-1}, showing that
\[\sum_{i=1}^{(p-1)/2} S_{2i}^* - S_{2i-1}^* \equiv - 1\pmod{p}.\]

\section{When do derivatives of distinct functions match?}
In this section, we explore the possibility for two functions  to have  their (higher-order) derivatives in different directions being functionally equal. That is we study functions the following problem.

\begin{prob}\label{prob:matching}
  Given two multisets $\{\alpha_1,\dots, \alpha_d\},\{\beta_1,\dots,\beta_d\}\subset \mathbb{F}_{p^n}^*$, what are the functions  $F,G:\mathbb{F}_{p^n}\rightarrow\mathbb{F}_{p^n}$ such that
  \begin{equation}\label{eq:ho-matching}
    \Delta_{\alpha_1,\dots,\alpha_d} F(x) = \Delta_{\beta_1,\dots,\beta_d} G(x), \quad \forall x\in\mathbb{F}_{p^n}?
  \end{equation}
\end{prob}

    We consider the two families of directions to have the same cardinality, since it is always possible to take $G(x)=\Delta_{\beta_{d+1}} G'(x)$. 
    
    Before we go further into the exploration of this problem, we recall the necessary and sufficient condition for a function to be a derivative in some direction $\alpha\in\mathbb{F}_{p^n}^*$:
    
    \begin{thm}[\cite{xiong2014} Theorem 1]
            Let $F:\mathbb{F}_{p^n}\rightarrow\mathbb{F}_{p^n}$. There exists a function $G:\mathbb{F}_{p^n}\rightarrow\mathbb{F}_{p^n}$ such that 
            \[ F(x) = \Delta_\alpha G(x)\]
            for some $\alpha\in\mathbb{F}_{p^n}^*$,
            if and only if
            \[ \Delta_{\underbrace{\alpha,\dots,\alpha}_{p-1}} F(x) = 0.\]
            
    \end{thm}
        
    And its natural corollary:
    \begin{cor}\label{cor:iff_cond}
        Let $F:\mathbb{F}_{p^n}\rightarrow\mathbb{F}_{p^n}$. There exists a function $G:\mathbb{F}_{p^n}\rightarrow\mathbb{F}_{p^n}$ such that 
            \[ F(x) = \Delta_{\underbrace{\alpha,\dots,\alpha}_{p-1}}  G(x)\]
            for some $\alpha\in\mathbb{F}_{p^n}^*$,
            if and only if
            \[  \Delta_\alpha F(x) = 0.\]
    \end{cor}

Two distinct functions having the same derivative in a given direction are related through Corollary~\ref{cor:iff_cond}:
\begin{prop}
Let $F$ and $G$ be two functions over $\mathbb{F}_{p^n}$ and $\alpha\in\mathbb{F}_{p^n}^*$. Then,
    \[
    \Delta_\alpha F(x) =  \Delta_\alpha G(x) \quad \Leftrightarrow \quad F(x) = G(x) + H(x)
    \]
    for any function $H:\mathbb{F}_{p^n}\rightarrow\mathbb{F}_{p^n}$ such that $\Delta_\alpha H(x) = 0$.
\end{prop}

Considering this last Proposition, we can now assume that the two families of directions in Equation~\ref{eq:ho-matching} are disjoint since otherwise we know how the functions $F$ and $G$ relate to each other. That is, we choose:
\[\{\alpha_1,\dots,\alpha_d\} \cap \{\beta_1,\dots,\beta_d\} = \emptyset. \]

We focus now on the first-order variant of the aforementioned Problem~\ref{prob:matching}, that is to find conditions on the functions $F,G:\mathbb{F}_{p^n}\rightarrow\mathbb{F}_{p^n}$ and directions $\alpha\neq\beta\in\mathbb{F}_{p^n}^*$ such that
\begin{equation}\label{eq:first-order}
  \Delta_{\alpha} F(x) = \Delta_\beta G(x).
\end{equation}
It is worth noting, that in this reduced variant of Problem~\ref{prob:matching}, there exists a link with antiderivative functions and so we recall a result from Salagean~\cite{Salagean20}.

\begin{thm}[\cite{Salagean20} Theorem 4]\label{thm:salagean}
  Let $\alpha_1,\dots,\alpha_k\in\mathbb{F}_{p^n}^*$ be $\mathbb{F}_p$-linearly independent elements and let $D_1,\dots,D_k:\mathbb{F}_{p^n}\rightarrow\mathbb{F}_{p^n}$. There exists a function $H:\mathbb{F}_{p^n}\rightarrow\mathbb{F}_{p^n}$ such that $\Delta_{\alpha_i} H(x) = D_i(x)$ for all $i=1,\dots,k$ if and only if both of the following conditions are satisfied:
  \begin{enumerate}
  \item $\Delta_{\underbrace{\alpha_i,\dots,\alpha_i}_{p-1}} D_i(x) = 0$
  \item $\Delta_{\alpha_j} D_i(x) = \Delta_{\alpha_i} D_j(x)$ for all $i,j=1,\dots,k$
  \end{enumerate}
\end{thm}

From this theorem, it is clear how the functions $F$ and $G$ from Equation~\eqref{eq:first-order} are related in the case where they are both already derivatives in the directions $\beta$ and $\alpha$ respectively. However, when the first condition of Theorem~\ref{thm:salagean} is not met, the relation between $F$ and $G$ is less obvious. We present here a simple construction satisfying only the second condition of the above Theorem.

In order to do so, we use the conversion between the univariate representation of a function $F:\mathbb{F}_{p^n}\rightarrow\mathbb{F}_{p^n}$ and its multivariate representation $F:\mathbb{F}_p^n\rightarrow\mathbb{F}_p^n$. We refer again to Salagean~\cite[Section 4.1]{Salagean20} for a detailed and comprehensive explanation of how it is done. In what follows, we make sure to use a basis $\mathcal{E}=\{e_1,\dots,e_n\}$ of $\mathbb{F}_p^n$ over $\mathbb{F}_p$ such that $\alpha,\beta$ are represented by canonical elements in the basis $\mathcal{E}$, say $e_1$ and $e_2$ respectively.

Let us define the concept of the degree in a variable of a multivariate function $f:\mathbb{F}_p^n\rightarrow\mathbb{F}_p$.
\begin{defn}
  Let a function $f:\mathbb{F}_p^n\rightarrow \mathbb{F}_p$ have an algebraic normal form
  \[ f(x_1,\dots,x_n) = \sum_{(i_1,\dots,i_n)\in\mathbb{F}_p^n} f_{(i_1,\dots,i_n)} \prod_{j=1}^n x_j^{i_j}.\]
  The degree in the variable $x_j$ of the function $f$ is:
  \[ \deg_{x_j}(f) = \max \{ i_j \mid f_{(i_1,\dots,i_n)} \neq 0,~ \forall (i_1,\dots,i_n)\in\mathbb{F}_p^n\}.\]
\end{defn}

\begin{rem}
  With the same notation as the above definition, note that
  \[ \deg_{x_i}(\Delta_{e_i} f) \leq \deg_{x_i}(f) - 1.\]
\end{rem}

Under such a multivariate representation, a vectorial function $F:\mathbb{F}_p^n\rightarrow\mathbb{F}_p^m$ is not a derivative function in a direction $e_i$ if and only if at least one of its component function has the maximal degree, $p-1$ in the variable $x_i$. Said differently and more formally:

\begin{prop}
  Let $F:\mathbb{F}_p^n\rightarrow\mathbb{F}_p^m$ and $\mathcal{E}=\{e_1,\dots,e_n\}$ be a basis of $\mathbb{F}_p^n$ over $\mathbb{F}_p$. The function $F$ is a derivative function in the direction $e_i$:
  \[\Delta_{\underbrace{e_i,\dots,e_i}_{p-1}} F(x_1,\dots,x_n) = 0,\]
  if and only if
  \[ \deg_{x_i}(\mu\cdot F) < p-1, \quad \forall \mu\in\mathbb{F}_p^m.\]
\end{prop}

We know that two vectorial functions are identical if and only if their component functions are pairwise identical. Hence, using the basis $\mathcal{E}$, and according to Equation~\eqref{eq:first-order}, if
\begin{equation}\label{eq:vectorial}
  \Delta_{e_1} F(x_1,\dots,x_n) = \Delta_{e_2} G(x_1,\dots,x_n),
\end{equation}
then clearly
\begin{equation}\label{eq:fo-component}
  \Delta_{e_1} \mu\cdot F(x_1,\dots,x_n) = \Delta_{e_2} \mu\cdot G(x_1,\dots,x_n),
\end{equation}
for any nonzero $\mu\in \mathbb{F}_p^n$.

We can deduce that the degree in $x_1$ of $\Delta_{e_1}f(x_1,\dots,x_n)$ cannot exceed $p-2$, and so is the degree in $x_2$ of $\Delta_{e_2} g(x_1,\dots,x_n)$. More generally,
\begin{lem}\label{lem:vardeg}
  Let $f,g:\mathbb{F}_p^n\rightarrow\mathbb{F}_p$ and $\mathcal{E}=\{e_1,\dots,e_n\}$ be a basis of $\mathbb{F}_p^n$ over $\mathbb{F}_p$ such that
  \[ \Delta_{e_i}f(x_1,\dots,x_n) = \Delta_{e_j} g(x_1,\dots,x_n),\]
  for some $i\neq j \in\{1,\dots,n\}$. Then
  \begin{align*}
    \deg_{x_j} (\Delta_{e_i} f) &= \deg_{x_j}(\Delta_{e_j}g) < p-1 \quad \mbox{ and} \\
    \deg_{x_i} (\Delta_{e_j} g) &= \deg_{x_i}(\Delta_{e_i}f) < p-1. \\
  \end{align*}
\end{lem}

We can finally state the main result of this section, consisting of a construction for pairs of functions satisfying the second condition of Theorem~\ref{thm:salagean} and answering the first-order variant of the Problem~\ref{prob:matching}.

    \begin{thm}\label{thm:fo-matching}
      Let $F,G:\mathbb{F}_p^n\rightarrow\mathbb{F}_p^n$ and $\mathcal{E}=\{e_1,\dots,e_n\}$ be a basis of $\mathbb{F}_p^n$ over $\mathbb{F}_p$. Then
      \[ \Delta_{e_1} F(x_1,\dots,x_n) = \Delta_{e_2} G(x_1,\dots,x_n), \quad \forall (x_1,\dots,x_n)\in\mathbb{F}_p^n\]

      if and only if for any nonzero $\mu \in\mathbb{F}_p^n$, the algebraic normal forms
      \begin{align*}
        \mu \cdot F(x_1,\dots,x_n) &= \sum_{(i_1,\dots,i_n)\in\mathbb{F}_p^n} f_{(i_1,\dots,i_n)}^{(\mu)} \prod_{j=1}^n x_j^{i_j} \quad \mbox{ and} \\
        \mu \cdot G(x_1,\dots,x_n) &= \sum_{(i_1,\dots,i_n)\in\mathbb{F}_p^n} g_{(i_1,\dots,i_n)}^{(\mu)} \prod_{j=1}^n x_j^{i_j}
        \end{align*}
      are such that
      \[ \sum_{k=i_1+1}^{p-1} \binom{k}{k-i_1} f_{(k,i_2,i_3,\dots,i_n)}^{(\mu)} = \sum_{k=i_2+1}^{p-1} \binom{k}{k-i_2} g_{(i_1,k,i_3,\dots,i_n)}^{(\mu)},\]
      for all $i_1,i_2=0,\dots,p-2$ and $(i_3,\dots,i_n)\in\mathbb{F}_p^{n-2}$, and
      \[ \sum_{k=j+1}^{p-1} \binom{k}{k-j} f_{(k,p-1,i_3,\dots,i_n)}^{(\mu)} = \sum_{k=j+1}^{p-1} \binom{k}{k-j} g_{(p-1,k,i_3,\dots,i_n)}^{(\mu)}=0 \]
      for all $j=0,\dots,p-2$ and $(i_3,\dots,i_n)\in\mathbb{F}_p^{n-2}$.

    \end{thm}
    \begin{proof}
      We have already established that Equation~\eqref{eq:vectorial} is satisfied if and only if Equation~\eqref{eq:fo-component} is satisfied for every $\mu\in\mathbb{F}_p^n$.
      
      We can re-write the algebraic normal forms of the component functions so that it emphasis the two variables we are interested in, namely $x_1$ and $x_2$.
      \begin{align*}
        \mu \cdot F(x_1,\dots,x_n) &= \sum_{i_1,i_2\in\mathbb{F}_p} x_1^{i_1}x_2^{i_2} \left( \sum_{(i_3,\dots,i_n)\in\mathbb{F}_p^{n-2}} f_{(i_1,\dots,i_n)}^{(\mu)} \prod_{j=3}^n x_j^{i_j}\right) \quad \mbox{ and} \\
        \mu \cdot G(x_1,\dots,x_n) &= \sum_{i_1,i_2\in\mathbb{F}_p} x_1^{i_1}x_2^{i_2} \left( \sum_{(i_3,\dots,i_n)\in\mathbb{F}_p^{n-2}} g_{(i_1,\dots,i_n)}^{(\mu)} \prod_{j=3}^n x_j^{i_j}\right).
      \end{align*}
      For the sake of simplicity we will denote $\bar{x}:=(x_1,\dots,x_n)$ and use the following notation for the part ne depending neither on $x_1$ nor $x_2$:
      \begin{align*}
        \sum_{w\in\mathbb{F}_p^{n-2}} f_{(i_1,i_2,w)}^{(\mu)} X_w &:= \sum_{(i_3,\dots,i_n)\in\mathbb{F}_p^{n-2}} f_{(i_1,\dots,i_n)}^{(\mu)} \prod_{j=3}^n x_j^{i_j} \quad \mbox{ and}\\
        \sum_{w\in\mathbb{F}_p^{n-2}} g_{(i_1,i_2,w)}^{(\mu)} X_w &:= \sum_{(i_3,\dots,i_n)\in\mathbb{F}_p^{n-2}} g_{(i_1,\dots,i_n)}^{(\mu)} \prod_{j=3}^n x_j^{i_j},
      \end{align*}
      and so the algebraic normal forms of the component functions become:
      \begin{align*}
        \mu \cdot F(\bar{x}) &= \sum_{i_1,i_2\in\mathbb{F}_p} x_1^{i_1}x_2^{i_2} \left( \sum_{w\in\mathbb{F}_p^{n-2}} f_{(i_1,i_2,w)}^{(\mu)} X_w\right) \quad \mbox{ and} \\
        \mu \cdot G(\bar{x}) &= \sum_{i_1,i_2\in\mathbb{F}_p} x_1^{i_1}x_2^{i_2} \left( \sum_{w\in\mathbb{F}_p^{n-2}} g_{(i_1,i_2,w)}^{(\mu)} X_w \right).
      \end{align*}
      
      The derivative in direction $e_1$ of $\mu\cdot F$ is 

      \begin{align*}
        \Delta_{e_1} \mu\cdot F(\bar{x}) &= \sum_{i_1\neq0,i_2\in\mathbb{F}_p} \left( (x_1+1)^{i_1} - x_1\right) x_2^{i_2} \left( \sum_{w\in\mathbb{F}_p^{n-2}} f_{(i_1,i_2,w)}^{(\mu)} X_w \right) \\
        &= \sum_{i_1\neq0,i_2\in\mathbb{F}_p} \left( \sum_{k=1}^{i_1} \binom{i_1}{k} x_1^{i_1-k} \right) x_2^{i_2} \left( \sum_{w\in\mathbb{F}_p^{n-2}} f_{(i_1,i_2,w)}^{(\mu)} X_w \right) \\
        &= \sum_{i_2\in\mathbb{F}_p} x_2^{i_2} \left( \sum_{w\in\mathbb{F}_p^{n-2}} f_{(1,i_2,w)}^{(\mu)} X_w + \left(\binom{2}{1}x_1+1\right)\sum_{w\in\mathbb{F}_p^{n-2}} f_{(1,i_2,w)}^{(\mu)} X_w  + \dots \right.\\
        &\left. \qquad + \left(\binom{p-1}{1}x_1^{p-2} + \binom{p-1}{2}x_1^{p-3} + \dots + 1 \right)\sum_{w\in\mathbb{F}_p^{n-2}} f_{(1,i_2,w)}^{(\mu)} X_w\right) \\
        &= \sum_{i_2\in\mathbb{F}_p} x_2^{i_2} \left(\sum_{w\in\mathbb{F}_p^{n-2}} X_w \left( f_{(1,i_2,w)}^{(\mu)} + f_{(2,i_2,w)}^{(\mu)} + \dots + f_{(p-1,i_2,w)}^{(\mu)}\right) \right.\\
        & \qquad + \sum_{w\in\mathbb{F}_p^{n-2}} x_1 X_w \left( \binom{2}{1}f_{(2,i_2,w)}^{(\mu)} + \binom{3}{2}f_{(3,i_2,w)}^{(\mu)} + \dots + \binom{p-1}{p-2}f_{(p-1,i_2,w)}^{(\mu)}\right) \\
        &\qquad + \sum_{w\in\mathbb{F}_p^{n-2}} x_1^2 X_w \left( \binom{3}{1}f_{(3,i_2,w)}^{(\mu)} + \binom{4}{2}f_{(4,i_2,w)}^{(\mu)} + \dots + \binom{p-1}{p-3}f_{(p-1,i_2,w)}^{(\mu)}\right) \\
        &\qquad + \qquad \dots \\
        &\left.\qquad + \sum_{w\in\mathbb{F}_p^{n-2}} x_1^{p-2} X_w \left( \binom{p-1}{1}f_{(p-1,i_2,w)}^{(\mu)}\right)\right)\\
        &= \sum_{j=0}^{p-2}\sum_{i_2\in\mathbb{F}_p} \sum_{w\in\mathbb{F}_p^{n-2}} x_1^jx_2^{i_2}X_w \left( \sum_{k=j+1}^{p-1} \binom{k}{k-j} f_{(k,i_2,w)}^{(\mu)}\right).
      \end{align*}

      Similarly,
      \begin{equation*}
        \Delta_{e_2} \mu\cdot G(\bar{x}) = \sum_{j=0}^{p-2}\sum_{i_1\in\mathbb{F}_p} \sum_{w\in\mathbb{F}_p^{n-2}} x_1^{i_1}x_2^{j}X_w \left( \sum_{k=j+1}^{p-1} \binom{k}{k-j} g_{(i_1,k,w)}^{(\mu)}\right).
      \end{equation*}

      From Lemma~\ref{lem:vardeg} we already know that the terms of highest degree in $x_1$ and $x_2$, that is $p-1$, vanish and so we deduce that:
      \[ \sum_{k=j+1}^{p-1} \binom{k}{k-j} f_{(k,p-1,w)}^{(\mu)} = \sum_{k=j+1}^{p-1} \binom{k}{k-j} g_{(p-1,k,w)}^{(\mu)}=0 \]
      for all $j=0,\dots,p-2$ and $w\in\mathbb{F}_p^{n-2}$.

      Now we can conclude that Equation~\eqref{eq:fo-component} is satisfied if and only if
      \[ \sum_{k=i_1+1}^{p-1} \binom{k}{k-i_1} f_{(k,i_2,i_3,\dots,i_n)}^{(\mu)} = \sum_{k=i_2+1}^{p-1} \binom{k}{k-i_2} g_{(i_1,k,i_3,\dots,i_n)}^{(\mu)},\]
      for all $i_1,i_2=0,\dots,p-2$ and $(i_3,\dots,i_n)\in\mathbb{F}_p^{n-2}$.
    \end{proof}

    Note that we can specify the condition that needs to be satisfied in order to find functions matching their first order derivatives without the existence of an antiderivative.

    \begin{cor}
      With the same notations as in Theorem~\ref{thm:fo-matching}, if there exists an element $\mu\in\mathbb{F}_p^n$ such that $f_{(0,p-2,i_3,\dots,i_n)}^{(\mu)} \neq 0$ for some $(i_3,\dots,i_n)\in\mathbb{F}_p^{n-2}$, then the function $F$ is not a derivative in $e_2$ and we can construct a function $G$ such that only the second condition of Theorem~\ref{thm:salagean} is satisfied. 
    \end{cor}

    Using Theorem~\ref{thm:fo-matching} we present an exemple of component functions having the same derivatives in directions $\mathbb{F}_p$-linearly independent but without the existence of an antiderivative.
    \begin{exmp}
      Let $p=3$ and $n=2$.
      We are looking for functions $f,g:\mathbb{F}_3^2\rightarrow\mathbb{F}_3$ such that
      \[\Delta_{(1,0)}f(x_1,x_2) = \Delta_{(0,1)} g(x_1,x_2)\] and
      $\Delta_{(0,1),(0,1)}f(x_1,x_2)$ is not the zero function.
      
      Let $f(x_1,x_2) = x_1^2x_2 + x_2^2$, then
      \begin{align*}
        \Delta_{(0,1),(0,1)} f(x_1,x_2) &= f(x_1,x_2) + f(x_1,x_2+1) + f(x_1,x_2+2) \\
        &= x_1^2x_2 + x_2^2 \quad + x_1^2(x_2+1) + (x_2+1)^2 \quad + x_1^2(x_2+2) + (x_2+2)^2 \\
        &= 2
      \end{align*}
      and
      \begin{align*}
        \Delta_{(1,0)}f(x_1,x_2) &= x_1^2(x_2+1) + (x_2+1)^2 \quad + 2x_1^2x_2 + 2x_2^2 \\
        &= 2x_1x_2 + x_2 \\
      \end{align*}
      Now by using Theorem~\ref{thm:fo-matching}, we have to solve the linear system corresponding to the coefficients of the equation $\Delta_{(1,0)}f(x_1,x_2) = \Delta_{(0,1)} g(x_1,x_2) $, that is:
      \begin{align*}
        0 &= (2g_{(0,1)} + g_{(0,2)}) \\
        0x_1 &= (2g_{(1,1)} + g_{(1,2)})x_1 \\
        0x_1^2 &= (2g_{(2,1)} + g_{(2,2)})x_1^2 \\
        1x_2 &= 2g_{(0,2)}x_2 \\
        2x_1x_2 &= 2g_{(1,2)}x_1x_2 \\
        0x_1^2x_2 &= 2g_{(2,2)}x_1^2x_2.
      \end{align*}

      In the end we find $g(x_1,x_2) = x_2 + 2x_1x_2 + 2x_2^2 + x_1x_2^2$.
    \end{exmp}

\section{Conclusion}

In this work, we propose a solution to a counting problem about the sum of elements in a subset of a prime field. A natural question arises concerning a similar analysis over different finite algebraic structures such as quotient rings $\mathbb{Z}/n\mathbb{Z}$, although the author is currently unaware of concrete applications.

\begin{open}
  Given a finite algebraic structure $\mathcal{A}$, what is the cardinality of the set
  \[ \{ S\subseteq \mathcal{A} \mid \# S = \ell, \quad \sum_{s\in S} s = u\}\]
  for some $0\leq \ell\leq \#\mathcal{A}$ and $u\in\mathcal{A}$?
\end{open}

Solving this problem over prime fields allows us to demonstrate an equivalent reprensentation for generalized differentials. We therefore start the study of distinct functions having the same derivative function but in different directions. Although we only cover the first-order derivative case, we give an explicit construction for such pair of functions. A natural way to extend these results is to find similar necessary and sufficient conditions for higher-order derivatives, starting with second-order derivatives.

\bibliographystyle{plain}
\bibliography{lset}

\end{document}